\input amstex
\documentstyle{amsppt} \magnification=1000
\hcorrection{-5mm}   \vcorrection{-6mm}   \hsize=16.5cm   \vsize=25cm
\pageno=1 \NoBlackBoxes \nologo \nopagenumbers

     \def\bur#1#2{ { \underset{#2} \to {#1} }}
\def\trm#1{ {\text{\rm{#1}}} } 
\def\ra{\rightarrow} \def\lra{\longrightarrow} \def\ms{\mapsto}
 
\def\dim{{\operatorname{dim}\,}}


\def\smooth{{\operatorname{smooth}}} \def\sing{{\operatorname{sing}}}
 
\def\hom{{\operatorname{hom}}}    
    
 \def\alg{{\operatorname{alg}}} 
 \def\supp{{\operatorname{support}}}
\def\cl{{c\ell}}

\def\AA{{\Cal A}}    \def\EE{{\Cal E}} \def\FF{{\Cal F}} \def\GG{{\Cal G}}
 \def\II{{\Cal I}}  \def\LL{{\Cal L}}  \def\OO{{\Cal O}} \def\RR{{\Cal R}}
 \def\VV{{\Cal V}} \def\WW{{\Cal W}}  \def\ZZ{{\Cal Z}} 
\def\N{\Bbb N} \def\Z{\Bbb Z}  \def\Q{\Bbb Q} \def\C{\Bbb C} \def\H{\Bbb H}  \def\P{\Bbb P}

\def\pt{\hskip0.5pt}

\def\prg#1#2{\vskip4mm\noindent{{\bf{\S \ #1.\qquad #2.}}}\vskip4mm}

\def\bproc#1{\noindent \bf #1.\hskip3mm\it}                              \def\eproc{\rm}
\def\proof{\parindent=0mm{\it\underbar{Proof}. \ } \parskip=1mm}         \def\qed{ \vskip-3pt {\rightline{$\square$}} \parskip=3mm \parindent=4mm }

\def\Thm#1{\noindent {\bf Theorem #1. } \it}
\def\Lemma#1{\noindent {\bf Lemma #1. } \it}
\def\rmk#1{\noindent{\bf Remark #1. }}

 at 8truept
 at 6truept
 at 5truept
 at 8truept
 at 8truept
 at 12truept

\def\nrs{1}  \def\gnt{2}  \def\gct{3}  \def\nrg{4}  \def\svt{5}
\def\rksdue{6}  \def\ntt{7} \def\lmm{8} \def\clm{9} \def\dga{10} \def\ztd{11} \def\ths{12} \def\rks{13} \def\thsc{14} \def\fct{15}
\def\nrc{16} \def\thu{17} \def\dgc{18} \def\prv{19} \def\clc{20} \def\rkss{21} \def\exm{22} \def\coro{23}

\topmatter

\title
Algebraic versus homological equivalence for singular varieties
\endtitle

\author
Vincenzo Di Gennaro, \ Davide Franco, \ Giambattista Marini
\endauthor

\leftheadtext{Vincenzo Di Gennaro, \ Davide Franco, \ Giambattista Marini}
\rightheadtext{Algebraic versus homological equivalence for singular varieties}

\abstract
Let $ \, Y \, \subseteq \, \P^N \, $ be a possibly singular projective variety, defined over the field of complex
numbers. Let $X$ be the intersection of $Y$ with $h$ general hypersurfaces of sufficiently large degrees.
Let $d>0$ be an integer, and assume that $\dim Y=n+h$ and $ \, \dim Y_\sing \, \le \, \min\{ d\!+\!h\!-\!1 , \, n\!-\!1 \pt \} \pt $. Let $Z$ be an
algebraic cycle on $Y$ of dimension $d+h$, whose homology class in $H_{2(d+h)}(Y;\,\Q)$ is non-zero. In the present paper we prove
that the restriction of $Z$ to $X$ is not algebraically equivalent to zero.
This is a generalization to the singular case of a result due to Nori in the case $Y$ is smooth.
As an application we provide explicit examples of singular varieties for which homological equivalence is different from the algebraic one.


\endabstract
\endtopmatter

\document
\baselineskip 12pt
\parskip=3mm   \parindent=4mm

\prg{0}{Introduction}

In 1910, in order to prove the so called \lq\lq fundamental
theorem of irregular surfaces\rq\rq, Poincar\'e introduced its
famous {\it normal functions} [P], [Mu, p. 9], [C, p. 73]. Using
which and Noether-Lefschetz Theory, in 1969 Griffiths provided the
first examples of smooth projective varieties  having
homologically trivial cycles that are not algebraically equivalent
to zero. They follow from the following Theorem [G]. {\it {Let $Y$
be a smooth complex projective variety of dimension $n+1=2r\geq
4$, $Z\neq 0\in H_{2r}(Y;\Q)$ be an algebraic cycle on $Y$ of
codimension $r$, and  $X$ be a general hypersurface section of $Y$
of degree $\gg 0$. Then $Z\cap X$ is not algebraically trivial on
$X$.}} A very striking result if we think that in 1958 Severi
believed that {\it "$\dots$ l'equivalenza topologica; la quale
per\`o forse coincide con l'equivalenza algebrica, a meno di
divisori dello zero. Cos\`\i \ accade di fatto per le curve di una
superficie"} [S, p. 46]. Next other special examples were found by
Ceresa [Ce], and Ceresa and Collino [CC]. A deep improvement of
Griffiths' result was made in 1993 by Nori [N, Theorem 1]. With a
new argument, based on his celebrated Connectivity Theorem, he
proved that on a smooth projective variety homological and
algebraic equivalence do not coincide in a wider range, where
Griffiths' methods do not work. All these results concern smooth
varieties, but the notion of algebraic equivalence, and more
generally a good part of the theory of algebraic cycles, can be
formulated also for singular varieties [Fu]. So it seems to us of
some interest to investigate such questions also in the singular
case.

Recently, by a suitable modification of Griffiths' argument,
relying on Intersection Homology Theory and on a monodromy result
for singular varieties [DF1] (see also [DF2], [DF3]), it has been
shown that Griffiths' Theorem holds true also when the ambient
variety $Y$ has  isolated singularities [DFM]. Continuing in this
direction, in the present paper we prove the following
generalization of the quoted theorem [N, Theorem 1] to the case
where $Y$ and its general hypersurface sections are singular.

\bproc{Theorem \nrs} Let $ \, Y \, \subseteq \, \P^N \, $ be a
projective variety (possibly singular) of dimension $ \, n+h \, $,
defined over the field of complex numbers,
and let $ \, d >0\, $ be an integer. Assume $ \, \dim Y_\sing \, \le \,
\min\{ d\!+\!h\!-\!1 , \, n\!-\!1 \pt \} \pt $. Then there exists
an integer $ \, c \, $ such that for  any $ \ a_i \, \ge \, c \pt $,
general sections $ \ X_i \, \in \, \vert \pt \OO_Y(a_i) \pt \vert \pt $, and for any
algebraic cycle $ \ Z \, \in \, \ZZ_{d+h}(Y) \ $ with non-trivial homology class
$ \ \cl(Z) \, \in \, H_{2(d+h)}(Y;\,\Q) \pt $, the restriction of the class
$[Z]\in CH_{d+h}(Y) $ to $ \dsize \ X \, := \, \bigcap_{i=1}^h \, X_i \, $ in
$CH_{d}(X)$ is not algebraically equivalent to zero.
\eproc

As an application of the theorem above we provide a class of examples of singular varieties where
homological and algebraic equivalence do not coincide (cfr.\;section \S\,3).

Our argument uses a desingularization $ \, f: \widetilde Y \lra Y \, $ of $ \, Y \pt $.
This leads us to consider a lifting $ \, \widetilde Z \, \in \, \ZZ_{d+h}(\widetilde Y) \, $ of $ \, Z \pt $,
namely a cycle satisfying $ \, f_*(\widetilde Z) \, = \, Z \pt $.
The idea of the proof consists of a analysis of the consequences of [N, Theorem 1] for the desingularization while taking $ \, f_* \pt $.
Roughly speaking, on one hand such a cycle $ \, \widetilde Z \, $ cannot be algebraically supported on the exceptional divisor $ \, E \pt $,
on the other hand, arguing by contradiction, up to conclude that it is not algebraically trivial down on $ \, X \cap Y_\smooth \pt $,
can be universally assembled as a cycle supported on $ \, E \pt $.

Aside with the above technical point, as we work with the pair $ \, (\widetilde Y , \, E) \, $, we need to prove
two slight generalizations of Nori's Theorem itself in the smooth case.
As the reader can image one of them is essentially a pair version (Theorem \nrc), the other one
concerns the difficulty that one needs to work with line bundles such as
$ \, f^*(\OO_Y(a)) \otimes \OO_{\widetilde Y}(D) \, $ (with $ \, D \pt $ fixed, supported on the exceptional divisor),
which are not exactly the same as sufficiently ample line bundles on $ \, \widetilde Y \, $ (see Theorem \gnt \,\, below).

The paper is organized as follows. In section \S\,1 we prove Theorem \gnt; in section \S\,2 we prove Theorem \nrs;
in section \S\;3, as an application of  Theorem \nrs,
we provide explicit examples of singular varieties for which homological equivalence is different from the algebraic one.

\prg{1}{\qquad A remark on Nori's Theorem}

In this section we deal with a technical point we shall meet in
the body of proof of our Theorem {\nrs} for singular varieties: we want to achieve a statement
(cfr. Theorem \gnt \ below) which is a slight generalization of [N, Theorem 1 and Theorem 2 (a)] in the smooth case.
Roughly speaking, such statement says that the constant appearing in the
Connectivity Theorem [N, Theorem 4], works unchanged also if we
tensor each $ \, \OO_Y(a_i) \, $ (cfr.\;Theorem 1) with a globally generated invertible sheaf $ \, \GG_i $.

\Thm{\gnt}
Let $ \, Y \, \subseteq \, \P^N \, $ be a smooth projective variety of dimension $ \, n+h \, $ defined over the field of complex numbers.
Let $ \ \GG_1, \, \dots, \GG_h \ $ be invertible sheaves on $ \, Y \, $ generated by global sections.
Then there exists an integer $ \, c \, $ such that the following holds:
for $ \ d \, > \, r \, \ge \, 0 \, $, $ \ a_1, \, \dots , \, a_h \, \geq \, c \pt $, $ \ g_i \, \geq \, 0 \pt $,
general sections $ \ X_i \, \in \, \vert \pt \OO_Y(a_i) \otimes \GG_i^{\otimes g_i}\pt \vert \, $,
and for any
algebraic cycle $ \ Z \, \in \, \ZZ_{d+h}(Y) \ $ with non-trivial homology class
$ \ \cl(Z) \, \in \, H_{2(d+h)}(Y;\,\Q) \pt $, the restriction of the class
$[Z]\in CH_{d+h}(Y) $ to $ \dsize \ X \, := \, \bigcap_{i=1}^h \, X_i \, $ in
$CH_{d}(X)$, does not belong to $ A_rCH_d(X) $.
\eproc

We refer to [N, p. 364] for the definition of the group $ A_rCH_d(X) $. Here we only
recall that  $ A_0CH_d(X)= CH_d(X)_\alg $.

Using the same argument as in [N, p. 366] (compare also with [Gre,
pp. 79-92], and [V, pp. 238-240]), previous theorem follows from a
corresponding adapted version of the quoted Connectivity Theorem.
Following [N], [Gre] and [V], put $S=\prod_{i=1}^{h}
\P(H^0(Y,\OO_Y(a_i)\otimes \GG^{\otimes g_i}))$, $A=Y\times S$,
$B=\{(y,f_1,\dots,f_h)\,\vert \,f_i(y)=0 \,\,\forall
i\,\,\}\subset A$, and $A_T=A\times _{S}T$ and $B_T=B\times _{S}T$
for any base change $T\to S$.

\Thm{\gct} With the above notation, for any number $c$,
there is a number $N(c)$ such that if $a_i\geq N(c)$ for all $i$,
then $H^k(A_T,B_T;\Q)=0$ for all $k\leq 2n$ and any smooth
morphism $T\to S$. \eproc

Since Mixed Hodge Theory involved in the Nori's proof of
Connectivity Theorem applies with no change in our more general
context, the proof of previous theorem reduces to the following
vanishing (compare with [Gre, loc. cit.]). Fix ample invertible
sheaves $\OO_Y(A_1),\dots, \OO_Y(A_h)$, and globally generated
invertible sheaves $\OO_Y(G_1),\dots, \OO_Y(G_h)$. Given integers
$a_1,\dots,a_h\geq 0$, $g_1,\dots,g_h\geq 0$, set
$$
E:=\oplus_{i=1}^{h} \OO_Y(a_iA_i+g_iG_i),
$$
and let $M_E$ be defined by the exact sequence $0\to M_E\to
H^0(Y,E)\otimes \OO_Y\to E\to 0$.

\Thm{\nrg}
For any integer $ \, c \ge 0 \, $ and any coherent sheaf $ \, \FF \, $ on $ \, Y $,
there is a constant $ \, n_0 = n_0(\FF, \, c) \, $ such that if $ \, a_i \ge n_0 \, $ for all $ i $, then
for any $ \, g_i \geq 0 , \ a > 0 , \ b > 0 \, $ one has $ \, H^a(Y, \FF\otimes E^{\otimes b}\otimes M_E^{\otimes c}) = 0 $.
\rm

In order to prove this, we need the following slight generalized version of the well-known Serre
Vanishing Theorem [H, p.228].

\Lemma{\svt}
Let $ \, Y \, $ be a smooth projective variety, $ \, \AA \, $ an ample invertible sheaf, $ \, \FF \, $ a coherent sheaf and
$ \, \GG \, $ be an invertible sheaf generated by global sctions.
Then there exists an integer $ \, n_0 \, $ not depending on $ \, \GG \, $ such that
$$
H^i (Y, \, \FF \otimes \AA^n \otimes \GG ) \quad = \quad 0 \ ,
\qquad \forall \ i \, \ge \, 0 \, , \ n \, \ge \, n_0.
$$
\rm

\proof There exists $ \, n_0' \, $, not depending on $ \, \GG \,
$, such that $ \, \LL \, = \, \AA^n \otimes \GG \otimes
\omega_Y^{-1} \, $ is ample for $n\geq n_0'$ (compare with [H, p.
169, Exercise 7.5, (a) and (b)]). Therefore $ \, H^i(Y, \, \AA^n
\otimes \GG ) \, = \, H^i(Y, \, \LL \otimes \omega_Y) \, = \, 0 \,
$ for $ \, i
> 0 \, , \ n \ge n_0' $, by Kodaira Vanishing Theorem [H, p.248]
(here we need $Y$ smooth). We now fix an embedding $ \, Y
\subseteq \P^N \pt $. As $ \, \FF \, $ is coherent there exists a
resolution as follows
$$
0 \, \lra \, \RR \, \lra \, \oplus \OO_Y(d_l) \, \lra \, \FF \, \lra \, 0
$$
with $ \, \RR \, $ coherent. Tensoring with $ \, \AA^n \otimes \GG \, $ we get the exact sequence
$$
\oplus H^i(Y, \OO_Y(d_l) \otimes \AA^n \otimes \GG) \, \lra \, H^i(Y, \, \FF \otimes \AA^n \otimes \GG) \, \lra \, H^{i+1}(Y, \, \RR \otimes \AA^n \otimes \GG).
$$
Then we are done because the left term vanishes for $ \, n \, \ge \, n_0'' \ $ (not depending on $ \, \GG $),
and the right term vanishes for $ \, n \, \ge \, n_0(\RR,\, i+1) \, $ (here we use descending induction on $ \, i \pt $).

\qed

We are in position to prove Theorem {\nrg}, hence Theorem {\gct} and Theorem {\gnt}.

\parindent=0mm{\it\underbar{Proof} \ $($of Theorem \nrg$)$.} \parskip=1mm
We follow the same argument as in [Gre, pp 79, 80]. Denote by $ \,
p_0, \, p_1, \dots, \, p_c \, $ the projections $ \, Y^{c+1} \ra Y
\, $ and set
$$
\EE:=p^*_1(E)\otimes\dots\otimes p^*_c(E).
$$
As in [Gre] one may construct a coherent sheaf $Z^0$ on $Y^{c+1}$,
depending only on $c$, and flat over $Y$ via $p_0$, such that
$$
p_{0_*}(Z^0\otimes \EE)=M_E^{\otimes c}. \tag *
$$
Hence also
$$
p_0^*(\FF\otimes E^{\otimes b})\otimes Z^0\otimes \EE
$$
is a coherent sheaf on $Y^{c+1}$, flat over $Y$ via $p_0$.  When
restricting it on the fibers $j_y:Y^c\subset Y^{c+1}$ of $p_0$, we
get a direct sum of sheaves like
$$
\left[j_y^*(p_0^*(\FF)\otimes Z^0)\right]\otimes
\left[q^*_1(\OO_Y(a_{i_1}A_{i_1}+g_{i_1}G_{i_1}))\otimes\dots \otimes q^*_c(\OO_Y(a_{i_c}A_{i_c}+g_{i_c}G_{i_c}))\right]
$$
with $ \, i_1,\dots,i_c\in\{1,\dots,h\} \ $(where $ \, q_j \, $ denotes the $ j $-projection $ \, Y^c \to Y $).
Applying Lemma {\svt} to the variety $ \, Y^c \, $ we deduce that the cohomology of these sheaves vanishes, and therefore
$$
R^ap_{0_*}\left(p_0^*(\FF\otimes E^{\otimes b})\otimes Z^0\otimes
\EE\right)=0
$$
for any $a>0$ [H, p. 290, Thm. 12.11]. By [H, p.252, Ex. 8.1], the
projection formula [H, p. 123-124, Ex. 5.1, (d)], and $(^*)$, we
deduce
$$
H^a(Y^{c+1}, p_0^*(\FF\otimes E^{\otimes b})\otimes Z^0\otimes
\EE)=H^a(Y, p_{0_*}\left(p_0^*(\FF\otimes E^{\otimes b})\otimes
Z^0\otimes \EE\right))=H^a(Y, \FF\otimes E^{\otimes b }\otimes
M_E^{\otimes c}).
$$
This concludes the proof of Theorem {\nrg} because, again by previous Lemma {\svt},
we know that
\vskip2mm
\centerline{$ H^a(Y^{c+1}, \, p_0^*(\FF\otimes E^{\otimes b}) \otimes Z^0 \otimes \EE) \, = \, 0 \ $.}
\qed

\prg{2}{The proof of  Theorem {\nrs}}

We start with a first easy reduction.

\rmk{\rksdue} It suffices to consider the case where $ \, Z \, $ is chosen a priori.
In fact, on one hand our constant $ \, c \, $ shall not depend on $ \, Z \, $,
on the other hand there are countable many algebraic families of cycles (so the general $ \, X \, $ works for all of them).

\medskip \noindent
{\bf Notation \ntt.} For the rest of
this section, $ \ Y \, \subseteq \, \P^N \ $ and $ \, Z \, $
denote respectively a (possibly singular) projective variety and a
cycle as in the statement of Theorem {\nrs}. Denote by ${Z_i}^+$
the irreducible components of $Z$ which appear in $Z$ with
positive multiplicity ${n_i}^+>0$, and by ${Z_i}^-$ the
irreducible components of $Z$ which appear with negative
multiplicity $-\,{n_i}^-<0$, so that $Z^+:= \sum {n_i}^+{Z_i}^+$
and $Z^-:= \sum {n_i}^-{Z_i}^-$ are effective, and $Z=Z^+-Z^-$.
For a sequence of integers $ \, a_1, \, ...,\, a_h \, $ we
consider the cartesian product $ \, \times \vert \pt \OO_Y(a_i)
\pt \vert \, $ and the corresponding family $ \, \vert \, X \,
\vert \, $ of intersections $ \ X \, = \, \bigcap \, X_i \pt $,
where $ \, X_i \, \in \, \vert \pt \OO_Y(a_i) \pt \vert \pt $. We
let $ \, \Omega \, \subseteq \, \times \P(H^0(Y,\, \OO_Y(a_i))) \, $ be the open subset of parameters corresponding to
$ \, n $-dimensional complete intersections $ \, X \pt $ in $ \, Y \, $ (in particular, for such $X$'s, the inclusion $X\subset Y$ is
a regular imbedding of codimension $h$ ([Fu], p. 437)), which
intersect properly each ${Z_i}^{\pm}$ in a subvariety $X\cap
{Z_i}^{\pm}$ of dimension $d$. By [Fu, Chapters 2 and 7] we know
that, for any $ \, t \, \in \, \Omega \pt  $,  the restriction of
the class $[Z]\in CH_{d+h}(Y) $ to $X_t$ in $CH_{d}(X_t)$ is
represented by a well defined  cycle, which we will denote by
$X_t\cap Z\in \ZZ_d(X_t)$.

Next lemma is fundamental for our purposes. It should be well known to experts,
nonetheless we did not succeed in finding an appropriate reference,
therefore we prove it (compare with [V2, Lemma 3.2]).

\Lemma{\lmm} The locus $ \Omega_0 $ of parameters $ t \in \Omega $
where the restriction of the class $ [Z] \in CH_{d+h}(Y) $ to $ X_t $
in $ CH_{d}(X_t) $ is algebraically trivial, is a countable union of
closed algebraic subsets of $ \Omega $. \rm

\proof Let $\ZZ_+(Y)$ be the Chow variety parametrizing effective
cycles of $Y$, and let ${\text{Hilb}}(\ZZ_+(Y))$ be the Hilbert
scheme of $\ZZ_+(Y)$ ([Ser], [K], [Fr], [E], [PK], [LF]). We will
adopt similar notation for $X_t$. For each $h\in
{\text{Hilb}}(\ZZ_+(Y))$, let $C_h\subseteq \ZZ_+(Y)$ be the
corresponding subscheme. Let $\Gamma_Y^0:=\{h\in
{\text{Hilb}}(\ZZ_+(Y))\,:\, h^0(C_h,\OO_{C_h})=1\}$. By
Semicontinuity Theorem [H, p. 288] we know that $\Gamma_Y^0$ is
open in ${\text{Hilb}}(\ZZ_+(Y))$. Denote by $\Gamma_Y$ the
closure of $\Gamma_Y^0$ in ${\text{Hilb}}(\ZZ_+(Y))$. By the
Principle of Connectedness [H, p. 281, Exercise 11.4] we deduce
that $C_h$ is connected, for any $h\in \Gamma_Y$. Let $\II$ be the
closed subset of ${\text{Hilb}}(\ZZ_+(Y))\times \ZZ_+(Y)\times
\Omega$ formed by those triples $(h,W,t)$ such that $h\in
\Gamma_Y\cap {\text{Hilb}}(\ZZ_+(X_t))$, $W\in\ZZ_+(X_t)$, and
such that $(X_t\cap Z^+)+W$ and $(X_t\cap Z^-)+W$ belong to
$C_h\subseteq \ZZ_+(X_t)$. Let $p$ be the projection
${\text{Hilb}}(\ZZ_+(Y))\times \ZZ_+(Y)\times \Omega\to\Omega$.
Since $p$ is projective, and $\II$ is closed, then also $p(\II)$ is closed.
As we have $ \, p(\II) = \Omega_0 \, $ [Fu, Example 10.3.3, p. 186], then we are done.
\qed

As a consequence of previous lemma, the proof of Theorem {\nrs} reduces to the following result.

\noindent
{\bf Claim \clm. } There exists $t\in\Omega$ such that the
restriction of $[Z]\in CH_{d+h}(Y) $ to $ \dsize \ X_t$ in
$CH_{d}(X_t)$ is not algebraically equivalent to zero.

We consider a desingularization
$$
f \ : \quad \tilde Y \ \lra \ Y
$$
and we let $ \ U \, := \, Y \smallsetminus Y_\sing \ $ be the smooth locus of $ \, Y \, $
and $ \, E \, \subseteq \, \tilde Y \, $ denote the exceptional divisor.
By abuse of notation we consider $ \, U \, $ also as the open subset $ \, \tilde Y \smallsetminus E \pt $.
There is a commutative diagram with exact rows ([Fu], p. 21):
$$
\CD \ \\
@. \ZZ_{d+h}(E) @> \tilde i_* >> \ZZ_{d+h}(\tilde Y) @> \tilde j^* >> \ZZ_{d+h}(U) @>>> \ 0 \ \\
@. @. @VV f_* V @VV = V \\
@. \ZZ_{d+h}(Y_\sing) @> i_* >> \ZZ_{d+h}(Y) @> j^* >> \ZZ_{d+h}(U) @>>> 0 \\ \ \\
\endCD
\tag{\dga}
$$
where $ \ j : \, U \hookrightarrow Y \, , \ i : \, Y_\sing \hookrightarrow Y \, ,
\ \tilde j : \, U \hookrightarrow \tilde Y \, , \ \tilde i : \, E \hookrightarrow \tilde Y \ $
denote the natural inclusions.

We note that $ \ \ZZ_{d+h}(Y_\sing) \, = \, 0 \ $ for dimensional
reason, so by the diagram one has that:
$$
\exists \ \ \tilde Z \ \ \in \ \ \ZZ_{d+h}(\tilde Y) \quad
\big\vert \quad f_*(\tilde Z) \ = \ Z \tag{\ztd}
$$
which is defined up to elements in $ \ \tilde i_* \big( \ZZ_{d+h}(E)\big) \pt $. We fix such a $ \, \tilde Z \pt $.

Next, the idea is that one to apply Theorem \nrc \ below. So, we premise a few considerations.

We let $ \, H \, $ be a hyperplane section of $ \, Y \subseteq \, \P^N \pt $. Then there exists a divisor $ \, D \, $
supported on the exceptional divisor $ \, E \, $ and $ \ a_0 \, \in \, \N \ $ such that
$$
f^*(a H) + D \ \trm{ is very ample on } \ \tilde Y \ \trm{ for all } \ a \, \ge \, a_0 \, .
$$
We fix such $ \, D \, $ and $ \, a_0 \pt $. Considering generic sections
$$
\tilde X_i \quad \in \quad \vert \, f^*(a_i H) + D_i \, \vert
\tag{\ths}
$$
with $ \ D_i = a_i' D \ $ (any $ a'_i > 0$) and $ \ a_i \, \ge \, a_0 a_i' \ $
(note the linear system above is still very ample)
and their images $ \ X_i \ = \ f \pt (\tilde X_i) \, $,
which we might note to be sections of the line bundle $ \, \OO_Y(a_i H) \, $
(in particular they are Cartier divisors), we put
$$
\tilde X \ = \ \bigcap_{i=1}^h \, \tilde X_i \quad ,  \qquad X \ = \ f(\tilde X) \ .
$$
As $ \, \tilde X \, $ is $ \, n $-dimensional and meets the exceptional divisor in dimension $ \, n-1 \, $
(hence the generic point of $ \, \tilde X \, $ is not contained in $ \, E $),
the generic point of $ \, X \, $ does not belong to $ \, Y_\sing \pt $.
Furthermore, since $ \, \dim Y_\sing \, \le \, n-1 \pt $, then one has the equality $ \ X \ = \ \bigcap_{i=1}^h \, X_i \, $
(in fact, $ \ \bigcap_{i=1}^h \, X_i \ $ cannot have components of dimension greater or equal than $ \, n \, $
contained in $ \, Y_\sing \, $).
Summing up, for a general $ \, \tilde X \, $ as above, we have that $ \, X \, := \, f(\tilde X) \, $ is parametrized for some $ \, t \in \Omega $.
As $ \, \tilde X \, $ is general, it meets $ \, \tilde Z \, $ transversally
and $ \, \tilde X \cap \tilde Z \, $ is a well-defined cycle in $ \ZZ_d(\tilde X) $,
which represents the restriction of the class $ \, [\tilde Z] \in CH_{d+h}(\tilde Y) \, $ in $ \, CH_{d}(\tilde X) $.
We also have $ \, f_*(\tilde X \cap \tilde Z) \, = \, X \cap Z $.

Defining $ \ U_{_X} \ = \ \tilde X \smallsetminus E \ = \ X \smallsetminus Y_\sing \, $, there is a diagram analog to diagram (\dga) ([Fu], 10.3.4, p. 186)
$$
\CD
@. CH_{d}(E \cap \tilde X)\diagup \alg @> \tilde i_* >> CH_{d}(\tilde X)\diagup \alg @> \tilde j^* >> CH_d(U_{_X})\diagup \alg @>>> \ 0 \ \\
@. @. @VV f_* V @VV = V \\
@. CH_{d}(Y_\sing \cap X)\diagup \alg @> i_* >> CH_{d}(X)\diagup
\alg @> j^* >> CH_{d}(U_{_X})\diagup \alg @>>> 0.
\endCD
$$
As we already remarked, since $ \ X \, = \, \cap X_i $, and the $ \, X_i \, $ are Cartier divisors,
one has a well-defined cycle $ \ X \cap Z \, \in \, \ZZ_{d}(X) \ $ satisfying $ \ f_*(\tilde X \cap \tilde Z) \, = \, X \cap Z \, $.

\rmk{\rks}
At this stage, excluding (by contradiction) that $ \, X \cap Z \not\in CH_{d}(X)_{\alg} $,
one would obtain $ \, j^*([X\cap Z]) \, = \, 0 \, $ and, in turn, from the diagram above one would infer that
$$
\tilde j^*\big([\tilde X \cap \tilde Z] \big) \quad \in \quad CH_d(U_{_X})_\alg.
\tag{\rks$'$}
$$
Therefore, there would exist $ \ W_{_X} \, \in \, \ZZ_{d}(E \cap \tilde X) \, $ such that
$ \, \tilde i_* \big(W_{_X} \big) \, \sim_{_{\alg}} \, \tilde X \cap \tilde Z  \, $ in $ \, \ZZ_d(\tilde X) $,
i.e. $ \ \tilde i_* \big([W_{_X}] \big) \, = \, [\tilde X \cap \tilde Z] \, \in \, CH_d(\tilde X)\diagup \alg \, $.
Then, as it will be clear from the \lq\lq End of proof" below, to deduce the existence of some $ \, X_t \, $ as in the Claim {\clm},
the problem would reduce, roughly speaking, to that one to lift the cycle $ \, W_{_X} \, $ to a cycle $ \, W \, \in \, CH_{d+h}(E)\otimes \Q \pt $.
This achievement, mutatis mutandi, shall be the core of the \lq\lq pair version" Theorem {\nrc} (stated and proved at the end of this section).

\proof{} (End of proof of Theorem {\nrs} using the \lq\lq pair version" Theorem {\nrc} below). \newline
\ In view of diagram (\dga), which passes through modulo Borel-Moore homology equivalence ([Fu, p. 371], [Fu2, p. 211]), one sees
that the cycle $ \, \tilde Z \, $ as in (\ztd) has non-trivial Borel-Moore homology class in $ \, H_{2(d+h)}^{^{B\!M}}(U;\,\Q) \pt $.
Thus we are in position to apply Theorem {\nrc} to our triple $ \ (\tilde Y , \, E , \, \tilde Z) \pt $. Namely
$$
\cl(Z) \, \bur{\ne}{{}_{(hp \ thm \ \nrs)}} \, 0 \ \in \
H_{2(d+h)}(Y;\,\Q)) \qquad =\!=\!\Longrightarrow \qquad \cl(\tilde
Z) \, \bur{\ne}{{}_{(hp \ thm \ \nrc)}} \, 0 \ \in \
H_{2(d+h)}^{^{B\!M}}(U;\,\Q).
$$
Applying Theorem {\nrc} one infers that there exists a constant $ \, c \, $ such that, for $ \, a_i' \ge c , \ d >  0 , \ g_i \ge 0 \, $ and general sections
$$
\tilde X_i \quad \in \quad \vert \, a_i' \big( f^*(a_0 H) + D
\big) + g_i f^*(H) \, \vert,
\tag{\thsc}
$$
the restriction of $ \, [\tilde Z] \, $ to $ \, U_{_X} \, $ is not algebraically equivalent to zero in  $ \, \ZZ_d(U_{_X})$.
Observe that the linear systems in (\thsc) are as in (\ths) for $ \, a_i \, = \, a'_i a_0 + g_i \, $ and $ \, D_i = a_i' D \, $ and, in turn,
for any $ \, a_i \, \ge \, a_0 c \pt $.
For constants $ \, a_i \, $ as above, arguing by contradiction, as already remarked, one
deduces the existence of some $ \, X_t \, $ (through the singular locus $ \, Y_\sing $, in fact our $ \, \tilde X \, $ cut $ \, E $) as in the Claim {\clm}.
Thus we are done with our reduction.
\qed

\rmk{\fct}
As observed $ \, X_i \, \in \, \vert a_i H \vert \pt $, using that $ \, [\tilde X\cap \tilde Z]\notin CH_{d}(\tilde X)_{\alg} \, $
we prove a Nori type result for the corresponding
$ \ X \, = \, \bigcap \pt X_i $ \ for any $ \, a_i \, \ge \, $ \lq\lq{\it fixed constant}" ($ = \, a_0 c \pt $).
We want to stress that without using Theorem {\nrc}, and in turn without using Theorem {\gnt}, one could only reach sequences of
integers (namely integers of type $ \, a_i' a_i \, $ with $ \, a_i' \ge c_i \, $ with no control of the $ \, c_i \pt $).

\bproc{Theorem \nrc \ (\lq\lq pair version of Nori's theorem")}
Let $ \, Y \, \subseteq \P^N \, $ be a smooth projective variety of dimension $ \, n+h \pt $,
and let $ \, E \subseteq Y \, $ be a simple normal crossing effective divisor.
Let $ \ U \, = \, Y \smallsetminus E \, $, and let  $ \, d > 0 \, $ be an integer.
Then there exists an integer $ \, c \, $ such that for  any $ \ a_i \, \ge \, c \pt $,
general sections $ \ X_i \, \in \, \vert \pt \OO_Y(a_i) \otimes \GG_i \pt \vert \, $
(with the $ \, \GG_i \, $ any globally generated line bundles),
and for any algebraic cycle $ \ Z \, \in \, \ZZ_{d+h}(Y) \ $ with non-trivial Borel-Moore homology class
$ \, \cl(Z) \, \in \, H_{2(d+h)}^{^{B\!M}}(U;\,\Q) \, $ (i.e. such that $ \, [Z]  \not\in i_* \, CH_{d+h}(E) + CH_{d+h}(Y)_\hom $),
one has that the restriction of $ \, [Z] \, $ to $ \, U_{_X} \, $ is not algebraically equivalent to zero in  $ \, \ZZ_d(U_{_X}) $,
where $ \dsize \ X \ = \ \bigcap_{i=1}^h \, X_i \pt $, $ \quad U_{_X} \, = \, X \smallsetminus E \pt $,
$ \, i: E \hookrightarrow \, Y $ and $ \, j : U_{_X} \hookrightarrow X \, $ denote inclusions.
\eproc

\proof
We construct a cycle $ \, W \, $ (or better, a multiple of it) as in the diagram below
thanks to the universal cycle associated with the various $ \, W_X \, $.

First we observe that we are free to fix a priori our $Z$ (compare with Remark {\rksdue}).
By Theorem {\gnt} we know that the (Noether-Lefschetz) general $ \, X \, $ satisfies the following property:
$$
\forall \ Z\,' \, \in \, \ZZ_{d+h}(Y) \ \big\vert \ Z\,' \, \not\sim_{_{\Q-\hom}} \, 0 \quad
\trm{ one has that } \quad \ [Z\,']\,\vert_ {_X} \, \ne \, 0 \ \in \ CH_d(X) \diagup \alg.
\tag{\thu}
$$
Consider the diagram
$$
\CD
\ \\
[W] @. \in @. CH_{d+h}(E)\diagup \alg @> i_* >> CH_{d+h}(Y)\diagup \alg @> j^* >> CH_{d+h}(U)\diagup \alg @>>> \ 0 \ @. \qquad \ \\
@. @. @VV _{ \cap X } V @VV _{ \cap X } V @VV _{ \cap X } V \\
@. @. CH_d(X\cap E)\diagup \alg @> i_* >> CH_d(X)\diagup \alg @> j^* >> CH_d (U_{_X})\diagup \alg @>>> 0 \\
[W_{_X}] @. \qquad @. @. \hskip-20mm \ms \hskip11mm  [Z ]\,\vert_ {_X} \hskip11mm \ms \hskip-20mm @. \bur{0}{\trm{(by contradiction)}} \\
\endCD
\tag{\dgc}
$$
and observe that, assuming by contradiction the vanishing $ \,
j^*( [Z ]\,\vert_ {_X}) \, = \, 0 \pt $ for some $Z$ and for a
general $X$, there exists a cycle $ \, W_{_X} \, $ as in the
diagram. In order to conclude, it suffices to lift it to a cycle $
\, W \pt \in \ZZ_{d+h}(E)\otimes \Q$. To see this, argue as
follows. The cycle $ \ Z' \ := \ Z - W \ $ satisfies $ \ \cl(Z')
\, \ne \, 0 \, \in \, H_{2(d+h)}(Y , \, \Q) \pt $, indeed $ \
\cl(Z') \vert_U \, = \, \cl(Z) \vert_U \, \ne \, 0 \pt $. On the
other hand from the diagram one obtains $[Z'] \vert_{_{X}} \quad
\in \quad CH_{d}(X)_\alg\otimes \Q$, therefore for some integer $m
> 0$ one has
$$
[mZ'] \vert_{_{X}} \quad \in \quad CH_{d}(X)_\alg \ .
$$
This contradicts (\thu).

Eventually, we come to the construction of our cycle $ \, W \pt $.
Put
$$
\eqalign{ \ZZ_d^{(E)} \quad & = \quad {\dsize
\bigcup_{\delta_1,\,\delta_2\, \in\, \N} } \, \VV_{\delta_1}
\times \VV_{\delta_2} },
$$
where $ \ \VV_\delta \ $ denotes the Chow variety of $ \, \delta $-degree and $ \, d \, $ dimensional effective cycles supported on
$ \, E \, $, and where to the point $ \ \sigma \, = \, (s_1,\,s_2) \, \in \, \VV_{\delta_1} \times \VV_{\delta_2} \ $ we associate
the cycle $ \ S_\sigma \, = \, S_{s_1} - S_{s_2} \, $, being $ \, S_{s_i} \, $ the cycle corresponding to $ \, s_i \in \VV_{\delta_i} $.
This exhibits $ \ZZ_d^{(E)} $ as a countable union of projective varieties parametrizing $d$-cycles on $E$. Denote by
$\H$  the Hilbert scheme parametrizing complete intersections in $Y$ of multidegree $(a_1,\dots,a_h)$, that meet $E$ and $Z$
properly, i.e. with the correct dimension (compare with Notation \ntt). For$ \ t \, \in \, \H \pt $, we denote by
$X_t$ the subscheme of $Y$ corresponding to the point $t$. Then we consider
$$
\Sigma_E \quad := \quad \big\{ \, (\sigma,\,t) \ \big\vert \
S_\sigma \, \in \, \ZZ_d (E\cap X_t) \, , \ i_* (S_\sigma) \,
\bur{\sim}{\bur{}{\bur{\trm{alg. equiv.}}{\trm{in } X_t}}} \, X_t
\cap Z \ \} \quad \subseteq \quad \ZZ_d^{(E)}\times \H.
$$
As in the proof of Lemma \lmm, we see that $\Sigma_E$ is a
countable union of quasi projective varieties, projective over
$\H$ via the natural projection $\pi_2:\Sigma_E\to \H$.
Taking into account that  we are arguing by contradiction, it
follows that $\pi_2:\Sigma_E\to \H$ is also surjective. As a
consequence one has that $ \ \Sigma_E \ $ contains a subvariety $
\ \Sigma\,'_E \ $ surjecting onto $\H$. We can fix a pair $ \,
(\sigma_{\!o}, \, o) \, $ as follows:
$$
(\sigma_{\!o}, \, o) \, \in \, \Sigma\,'_E \quad \big\vert \quad
X_o \ \trm{ is general (w.r.t. \thu) and meets $Z$ and  $E$
transversally.}
$$
Such $ \ S_{\sigma_{\!o}} \ $ shall be our $ \ W_{_{X_o}} \ $
(diagram \dgc). Let $\H^{(h)}$ be a sufficiently general
projective subvariety of $\H$ of dimension $h$, passing
through $o$ (e.g. take, for any $i=1,\dots,h$,   a general $
\P^1\subset \vert \pt \OO_Y(a_i) \pt \vert$, and define $\H^{(h)}$ as the image of $\times_{i=1}^{h} \P^1$ in $\H$), and
let ${\Sigma\,'_E}^{(h)}\subset \Sigma\,'_E$ be a $h$-dimensional
irreducible projective subvariety through the point
$(\sigma_{\!o}, \, o) \, \in \, \Sigma\,'_E$, mapping onto $\H^{(h)}$ via $\pi_2$
(so that $ \, \pi_2\vert_{{\Sigma\,'_E}^{(h)}} : \, {\Sigma\,'_E}^{(h)} \to \H^{(h)} \, $ is a generically finite map between projective varieties). Now consider
$$
\VV \quad := \quad \big\{ \, (p,\,\sigma,\,t) \ \big\vert \
(\sigma,\,t) \, \in \, {\Sigma\,'_E}^{(h)} \, , \ p \, \in \,
\supp (S_\sigma) \ \big\} \quad \subseteq \quad Y \times
{\Sigma\,'_E}^{(h)},
$$
and observe that such $ \, \VV \, $ is projective over $Y$, via
the projection $\pi_1: Y\times {\Sigma\,'_E}^{(h)}\to Y$.

\vskip3mm
{\bf Claim \prv. } For $ \ t \ \in \ \H^{(h)} \ $ the following properties hold:

{\quad \bf a) \ }
$ \dim \VV \ = \ d+h \, $;

{\quad \bf b) \ } $ \VV \vert_{Y \times \{(\sigma,\,t) \} } \quad
= \quad S_{\sigma}; $

{\quad \bf c) \ }
$ (\pi_1)_* \, \VV \ $ \ is supported on $ \ E \, $;

{\quad \bf d) \ }
$ \big[ (\pi_1)_* \, \VV \big] \big\vert_{X_t} \ =
\ {\dsize \bigcup_{(\sigma,\,t) \, \in \, {\Sigma\,'_E}^{(h)} \, \vert \ \pi_2(\sigma,\,t) \, = \, t } } \, S_\sigma \ $
(recall that $ \, \pi_2\vert_{{\Sigma\,'_E}^{(h)}} : \, {\Sigma\,'_E}^{(h)} \to \H^{(h)} $ is generically finite);

{\quad \bf e) \ }
for $ \ W \ := \ \dsize {1 \over m} \ (\pi_1)_* \, \VV \ $  and $ \ m \, = \,
\deg \left(\, \pi_2\vert_{{\Sigma\,'_E}^{(h)}} : \, {\Sigma\,'_E}^{(h)} \to \H^{(h)}\right)\pt $, \ one has
$$
W \cap X_t \, = \, \ W_{_{X_t}} \ .
$$

The crucial point is the one of proving $ d) $. To this purposes, it suffices to observe the following.
The family $ \, \{X_t\}_{t\in \H^{(h)}} \, $ contains the family $ \, \{W_{_{X_t}}\}_{t\in \H^{(h)}} \, $,
and it is cut out by $ \, X_o \, $ in $ W_{_{X_o}} $;
on the other hand, since $\H^{(h)}$ is sufficiently general, at a generic point in $ \, W_{_{X_o}} \pt $,
locally in $ \ \big(q,\, ( W_{_{X_o}} , \, o )\big) \ $ the family $ \, \{X_t\}_{t\in \H^{(h)}} \, $ is isomorphic
to $ \, Y \pt $, in other terms $ \, \pi_1 \, $ is locally at $ \ \big(q,\, ( W_{_{X_o}} , \, o )\big) \ $ an isomorphism.
Therefore, for $ \, t \ne o \, $, the variety $ \, X_t \pt $ meets $ \, W_{_{X_o}} \, $ transversally (in dimension $ \, d\!-\!h $),
and the local contribution of $ \ \VV \ $ in $ \ ( W_{_{X_o}} , \, o ) \, \in \, {\Sigma\,'_E}^{(h)} \ $ is $ \ W_{_{X_o}} \, $, i.e.
there exists an open neighborhood $ \ {\WW} \, \subseteq \, {\Sigma\,'_E}^{(h)} \ $ of the point $ (W_{_{X_o}}, \, o ) \, $ such that
$$
\pi_1\big( \VV \cap \big (Y \times {\WW}
\big) \ \big) \ \bigcap \ X_o \quad = \quad \ W_{_{X_o}}
\tag{\clc}
$$
and the intersection is transverse.
We may require that the statement (\clc) holds locally in $ \, o \, \in \, \H^{(h)} \, $
(via the generically finite morphism $ \, {\Sigma\,'_E}^{(h)} \, \ra \, \H^{(h)} \pt $),
rather than locally in $ \, (W_{_{X_o}}, \, o) \, \in \, {\Sigma\,'_E}^{(h)} \pt $.
This concludes our proof.
\qed

{\bf{Caution \rkss.}}
The previous considerations on the universal cycle $ \ \VV \ $ does not imply, as it could seem,
the possibility to lift any cycle on a general $ \ X_t \ $ to $ \ Y \ $
(namely the surjectivity of the restriction map \quad $ CH_{d+h}(Y)\slash\alg \otimes \Q \to CH_d(X_t)\slash\alg \otimes \Q $).
In our particular context, the restriction of $ \VV $  to $ X_t $ gives $ m $ cycles
algebraically equivalent to each other (or better to $ X_t \cap Z $),
but in general one only knows that such cycles are homologically equivalent [Mo].

\prg{3}{Example where algebraic and homological equivalence do not
coincide}

Fix a subspace $S$ of dimension $s$  in $\P:=\P^{n+h+1}$, with
$-1\leq s\leq \frac{n+h}{2}-1$. Let $T_1,\dots,T_q$ be general
subspaces of $\P^{n+h+1}$ of dimension $\frac{n+h}{2}$ containing
$S$, spanning $\P^{n+h+1}$ (a fortiori $q\geq 2$), and such that
$T_i\cap T_j=S$ when $i\neq j$. Let $Y\in
|H^0(\P^{n+h+1},\II_T(\ell))|$ be a general hypersurface in
$\P^{n+h+1}$ of degree $\ell\gg 0$ containing $T:=T_1\cup\dots
\cup T_q$. Observe that $T^0:=T\backslash S$ is smooth of
dimension $\frac{n+h}{2}$.

\noindent
{\bf{Proposition {\exm}}}.

1) $Y$ is irreducible, and $Y_\sing=S$.

2) $T_1,\dots,T_q$ are linearly independent in $H_{n+h}(Y,\Z)$. In particular $ \cl (Z:=T_1-T_2) \neq 0 \in H_{n+h}(Y,\Z)$.

3) For a general complete intersection $X\subset Y$ of dimension $n$, one has $ \cl(Z\cap X) = 0 \in H_{2d}(X,\Z)$ ($2d=n-h>0$).

{\proof} 1)
Since the tangent cone to $ \, S \, $ is the whole projective space $ \, \P $, it suffices to prove that $ Y_\sing \subseteq S $.
By Bertini's Theorem we know that $Y_\sing\subseteq T$.
Therefore we only have to check that  $Y$ is regular at any point of $ T^0 := T\backslash S $.
We will follow [OS, proof of Theorem 1.2] (compare also with [KA]).

Set $\Pi:=H^0(\P^{n+h+1},\II_T(\ell))$, and let $\Sigma\subseteq
\Pi\times T^0$ be the set of pairs $(f,x)$ such that $x\in
\sing(f)$. For any $x\in T^0$ set $\Pi_x:=\Sigma\cap (\Pi\times
\{x\})$. Recall that
$$
T^*_{x,T^0}=T^*_{x,T}={\goth m}_{x,\P}\slash \II_{x,T}+{\goth
m}_{x,\P}^2=\frac{{\goth m}_{x,\P}\slash{\goth m}_{x,\P}^2}{(\II_{x,T}+{\goth m}_{x,\P}^2)\slash {\goth m}_{x,\P}^2}.
$$
Since $\Pi$ generates ${(\II_{x,T}+{\goth m}_{x,\P}^2)\slash
{\goth m}_{x,\P}^2}$, and $\Pi_x$ maps to zero, we have a natural
surjection
$$
\Pi\slash \Pi_x\to {(\II_{x,T}+{\goth m}_{x,\P}^2)\slash {\goth
m}_{x,\P}^2}\to 0.
$$
It follows that
$$
\dim T^*_{x,T^0}=\frac{n+h}{2}=n+h+1-\dim ({(\II_{x,T}+{\goth
m}_{x,\P}^2)\slash {\goth m}_{x,\P}^2})\geq n+h+1-\dim
\Pi+\dim\Pi_x,
$$
hence
$$
\dim\Pi_x\leq \dim\Pi-\frac{n+h}{2}-1.
$$
Considering the projection $\Sigma\to T^0$, we deduce that
$$
\dim \Sigma\leq \dim T^0+\dim \Pi_x\leq \dim \Pi-1.
$$
Therefore the image of $\Sigma\to \Pi$ has dimension strictly less
than $\dim \Pi$.

2) First consider the exact sequence
$$
H_{n+h+1}^{BM}(Y\backslash T;\Z)\to H_{n+h}(T;\Z)\to H_{n+h}(Y;\Z)
$$
(recall that in the compact case Borel-Moore and singular homology
agree [Fu2, p. 217]). Since $H_{n+h}(T;\Z)$ is freely generated by
the components of $T$ [Fu2, p. 219, Lemma 4], it suffices to prove
$H_{n+h+1}^{BM}(Y\backslash T;\Z)=0$. To this purpose first
observe that $H_{n+h+1}^{BM}(Y\backslash
T;\Z)=H^{n+h-1}(Y\backslash T;\Z)$ because $Y\backslash T$ is
smooth [Fu2, p. 217]. Moreover by Lefschetz Theorem with
Singularities [GMP, p. 199], we know that
$H^{n+h-1}(\P^{n+h+1}\backslash T;\Z)\cong H^{n+h-1}(Y\backslash
T;\Z)$. It remains to check that $H^{n+h-1}(\P^{n+h+1}\backslash
T;\Z)=0$. This follows considering the exact sequence
$$
H_{n+h+3}(\P^{n+h+1};\Z)\to H^{BM}_{n+h+3}(\P^{n+h+1}\backslash
T;\Z)\to H_{n+h+2}(T;\Z),
$$
and observing that $H_{n+h+3}(\P^{n+h+1};\Z)=0$ ($n+h+3$ is odd),
as before one has $H^{n+h-1}(\P^{n+h+1}\backslash T;\Z)\cong
H^{BM}_{n+h+3}(\P^{n+h+1}\backslash T;\Z)$, and
$H_{n+h+2}(T;\Z)=0$ because $n+h+2>2\dim_{\C} T$ [Fu2, p. 219,
Lemma 4].

3) It is enough to prove that $H_{2d}(X;\Z)\cong H_{2d}(\P^{n+h+1};\Z)$ via push-forward.
This follows again by Lefschetz Theorem with Singularities [GMP, p. 199], because $X$ is a complete intersection of dimension $n\geq 2d+1$.
\qed

Here is an application of our main result Theorem 1.

{\bf{Corollary {\coro}}}. Let $X\subset Y\subset \P^{n+h+1}$ be as
above, with $s=\frac{n+h}{2}-1$. Then $Z\cap X$ is homologous to
zero in $X$, it is not algebraically zero, and $\dim X_\sing=d-1$.

\enddocument